\documentclass[12pt,noamsfonts]{amsart}

\usepackage{amsfonts}
\usepackage{amssymb}
\usepackage{amscd}

\newcommand{\marginlabel}[1]%
  {\mbox{}\marginpar{\raggedleft\hspace{0pt}\bfseries\sf#1}}
\def\NN{{\mathbb N}}

\def\CC{{\mathbb C}}
\def\AA{{\mathbb A}}
\def\RR{{\mathbb R}}
\def\QQ{{\mathbb Q}}

\def\cI{\mathcal{I}}

\def\cO{\mathcal{O}}

\DeclareMathOperator{\codim}{codim}

\DeclareMathOperator{\Hom}{Hom}

\DeclareMathOperator{\Spec}{Spec}

\DeclareMathOperator{\ord}{ord}
\DeclareMathOperator{\mld}{mld}

\newtheorem{lemma}{Lemma}[section]
\newtheorem{theorem}[lemma]{Theorem}
\newtheorem{corollary}[lemma]{Corollary}
\newtheorem{proposition}[lemma]{Proposition}

\theoremstyle{definition}
\newtheorem{definition}[lemma]{Definition}
\newtheorem{remark}[lemma]{Remark}
\newtheorem{conjecture}[lemma]{Conjecture}

\theoremstyle{remark}
\newtheorem*{remark*}{Remark}
\newtheorem*{note*}{Note}

\begin{document}

\title{Jet schemes, log discrepancies and inversion of adjunction}

\author[L. Ein]{Lawrence Ein}
\address{Department of Mathematics, Statistics and
Computer Science, University of Illinois at Chicago,
851 Morgan St., M/C. 249, Chicago, IL 60607-7045, USA}
\email{{\tt ein@math.uic.edu}}

\author[M. Musta\c{t}\v{a}]{Mircea~Musta\c{t}\v{a}}
\address{Department of Mathematics, Harvard University,
One Oxford Street, Cambridge, MA 02138, USA}
\email{{\tt mirceamustata@yahoo.com}}

\author[T. Yasuda]{Takehiko~Yasuda}
\address{Department of Mathematical Sciences, University of Tokyo, Komaba,
Meguro, Tokyo, 153-8914, Japan}
\email{{\tt t-yasuda@ms.u-tokyo.ac.jp}}

\thanks{The first author was partially
supported by NSF grant DMS 0200278.
The second author
served as a Clay Mathematics Institute Long-Term Prize Fellow
while this research has been done. The third author was financially
supported by the Japan Society for the Promotion of Science.}

\keywords{Jet schemes, motivic integration, minimal log discrepencies,
inversion of adjunction}

\maketitle

\section*{Introduction}

Singularities play a key role in the Minimal Model Program.
In this paper we show how some of the open problems in this area
can be approached using jet schemes.

Let $(X,Y)$ be a pair, where $X$ is a $\QQ$-Gorenstein normal 
variety, and $Y$ stands for a formal combination
$\sum_{i=1}^kq_i\cdot Y_i$, where $q_i\in\RR_+$
and $Y_i\subset X$ are proper closed subschemes. 
Fix a closed
subset $\emptyset\neq W\subseteq X$.
Using a suitable resolution of
singularities for the pair $(X,Y)$ one can define numerical invariants
$\mld(W;X,Y)$, called minimal log discrepancies. These invariants
in turn can be used to define the classes of singularities which appear
in Mori Theory.

We provide a way to compute minimal log discrepancies using arcs and jets. 
The $m$th jet scheme $X_m$ of $X$ is given set-theoretically as
$\Hom(\Spec\CC[t]/(t^{m+1}),X)$.
The limit of these schemes is the space of arcs 
$X_{\infty}=\Hom(\Spec\CC[[t]],X)$.
This is an infinite dimensional space, but we may associate to $(X,Y)$ a family
of subsets of $X_{\infty}$ of finite codimension. Given $W$, if we restrict
these subsets over $W$, then from their codimensions we can compute
$\mld(W;X,Y)$. We stress that this characterization holds in complete
generality. It extends the results in \cite{mustata1}, \cite{yasuda},
and \cite{elm},
where criteria were given for having non-negative log discrepancy,
under certain hypotheses on the singularities of $X$. The main ingredient 
in the proof of this characterization
is the theory of motivic integration on singular varieties, developed
by Denef and Loeser in \cite{denef}.

Note that our setting is slightly different than the 
standard one in Mori Theory. The usual setting is that of a pair
$(X,D)$, where $X$ is a normal variety, and $D$ is a $\QQ$-divisor
such that $K_X+D$ is $\QQ$-Cartier. We mention that our characterization
of minimal log discrepancies has an
analogue in this context (see Remark~\ref{classical}).
However, for the approach via spaces of arcs, our setting
seems more suggestive.

We leave the 
precise statement of our characterization for the main body of the 
paper (see Theorem~\ref{char_log_canonical}) 
and describe the consequences. Our first application is
a precise version of the
Inversion of Adjunction Conjecture
of Koll\'{a}r and Shokurov, in the case when the ambient variety is
smooth.

\begin{theorem}\label{inv_adj0}
Let $X$ be smooth, $Y=\sum_i q_i\cdot Y_i$ 
as above, and $D\subset X$
a normal effective  divisor such that 
$D\not\subseteq\bigcup_iY_i$.
For every proper closed subset $W\subset D$, we have
$$\mld(W;X,D+Y)=\mld(W;D,Y\vert_D).$$
\end{theorem}

Koll\'{a}r, Shokurov and Stevens 
proved special cases of Inversion of
Adjunction: see \cite{kollar}, \cite{shokurov1} and 
\cite{stevens}. The traditional approach 
to this problem involves applications of vanishing theorems.
We refer to \cite{kollar} for this part of the story.
We mention also the result of Ambro \cite{ambro1}
who proved Inversion of Adjunction in the case when
$X=\AA^n$ and $D$ is a hypersurface which is general
with respect to its Newton polyhedron.

 Together with our characterization
of minimal log discrepancies, Theorem~\ref{inv_adj0}
can be used to characterize terminal hypersurface singularities.
Recall that if $V$ is a locally complete intersection variety, it is proved in
\cite{mustata2} that $V$ has canonical singularities if and only if
$V_m$ is irreducible for all $m$. At least in the case when $V$
is a hypersurface in a smooth variety, 
this follows also from the above results.
Moreover, we get a similar characterization for the terminal case,
which was suggested by Mirel Caib\v{a}r.

\begin{theorem}\label{terminal0}
Let $X$ be smooth and $D\subset X$ an irreducible and reduced divisor.
$D$ has terminal singularities
if and only if $D_m$ is normal for every $m$.
\end{theorem}

Our final application is towards a semicontinuity statement.
Shokurov has given in \cite{shokurov} a conjectural uniform bound
for minimal log discrepancies. Ambro has made a stronger conjecture in
\cite{ambro} and he showed that this conjecture is equivalent to 
a semicontinuity statement about log discrepancies.
We prove this conjecture in the case of an ambient smooth variety.

\begin{theorem}
Let $X$ be a smooth variety, and $Y=\sum_iq_i\cdot Y_i$
as above. The function
$x\in X\longrightarrow\mld(x;X,Y)$ is lower semicontinuous.
\end{theorem}

\smallskip

A few words about the structure of the paper: in the first section we review 
the basic definitions and properties of minimal log discrepancies,
while in the second section we prove our characterization of these
invariants. In the next section we study the jet schemes of a hypersurface
in a smooth variety, and as a result, we prove Theorems~\ref{inv_adj0}
and \ref{terminal0}. In the last section we prove the above semicontinuity
statement.

\bigskip

\subsection{Acknowledgements}
We are  indebted to Florin Ambro,
Mirel Caib\v{a}r, Fran\c{c}ois Loeser, Rob Lazarsfeld, 
Mihnea Popa, and Vyacheslav V.~Shokurov for useful discussions and suggestions.
The third author thanks Yujiro Kawamata for his encouragement.
The paper has also benefitted from the referee's comments.
This work has been done while the second and the third author were visiting
Isaac Newton Institute for Mathematical Sciences; they are grateful for
hospitality.

\section{Log discrepancies}

All our varieties are defined over $\CC$. In this section
we review the definition and the basic properties of log discrepancies.
For a detailed discussion and proofs we refer to \cite{ambro}. 
Note that unlike in \cite{ambro}, in this paper we 
allow pairs of arbitrary codimension, but all the proofs can be reduced
to the case of divisors.

We will always work in the following setting. Let $X$ be a normal,
$\QQ$-Gorenstein variety, and $Y$ a formal combination
$Y=\sum_{i=1}^k q_i\cdot Y_i$, where $q_i\in\RR$, and where $Y_i\subset X$
are proper closed subschemes. We will restrict later to the case
when $q_i\geq 0$ for all $i$.
A divisor $E$ over $X$ is a prime Weil divisor
on $X'$, for some normal variety $X'$, proper and birational over $X$. We 
identify $E$ with the corresponding valuation of the function field of $X$.
The center of this valuation on $X$ is denoted by $c_X(E)$.

Given a divisor $E$ over $X$, choose a proper, birational morphism
$\pi : X'\longrightarrow X$, with $X'$ normal and $\QQ$-Gorenstein,
such that $E$ is a Cartier divisor on $X'$, and 
such that all the scheme-theoretic inverse images $\pi^{-1}(Y_i)$ 
are Cartier divisors. We write $\pi^{-1}(Y):=\sum_iq_i\cdot \pi^{-1}(Y_i)$.
The coefficient of $E$ in $K_{X'/X}-\pi^{-1}(Y)$
is $a(E;X,Y)-1$. It is clear that $a(E;X,Y)$ does not
depend on the particular model $X'$ we have chosen.

\begin{definition}\label{def1}
Let $W\subseteq X$ be a nonempty closed subset. The minimal log discrepancy
of $(X,Y)$ on $W$ is defined by
\begin{equation}
\mld(W;X,Y):=\inf_{c_X(E)\subseteq W}\{a(E;X,Y)\}.
\end{equation}
\end{definition}

\begin{definition}
The pair $(X,Y)$ is called log canonical if 
we have $\mld(X;X,Y)\geq 0$.
\end{definition}

We collect in the next proposition a few well-known facts about
minimal log discrepancies.

\begin{proposition}\label{basic_properties}
Let $(X,Y)$ and $W\subseteq X$ be as above.
\item{\rm (i)} If $W_1,\ldots,W_r$ are the irreducible components of $W$,
then $\mld(W;X,Y)=\min_i\mld(W_i;X,Y)$.
\item{\rm (ii)} If 
$U\subset X$ is an open subset such that both $W\cap U$ and $W\setminus U$
are nonempty, then 
$$\mld(W;X,Y)=\min\{\mld(W\setminus U;X,Y),
\mld(W\cap U;U,Y\vert_U)\}.$$
\item{\rm (iii)} $\mld(W;X,Y)\geq 0$ if and only if there is an open
subset $U$ of $X$, such that $W\subseteq U$ and $(U,Y\vert_U)$
is log canonical. If $\dim\,X\geq 2$, and if $\mld(W;X,Y)<0$,
then $\mld(W;X,Y)=-\infty$.
\item{\rm (iv)} If $\pi : X'\longrightarrow X$ 
is a proper, birational morphism,
with $X'$ normal and $\QQ$-Gorenstein, then 
$$\mld(W;X,Y)=\mld(\pi^{-1}(W);X',\pi^{-1}(Y)-K_{X'/X}),$$
where $\pi^{-1}(Y)=\sum_iq_i\cdot\pi^{-1}(Y_i)$.
\end{proposition}

The next proposition shows that minimal log discrepancies can be computed
using log resolutions. Given $(X,Y)$ and $W\subseteq X$, consider
$\pi : X'\longrightarrow X$ proper, birational, with $X'$ smooth, such that
$\pi^{-1}(Y)\cup {\rm Ex}(\pi)$ is a divisor with
 simple normal crossings. Here ${\rm Ex}(\pi)$ denotes the
exceptional locus of $\pi$. In addition, if
$W\neq X$, we assume that $\pi^{-1}(W)\cup\pi^{-1}(Y)\cup {\rm Ex}(\pi)$
is also a divisor with simple normal 
crossings. Note that by \cite{hironaka}, we can always find such a morphism.
Let us write $K_{X'/X}-\pi^{-1}(Y)=\sum_i(a_i-1)E_i$.

\begin{proposition}\label{computation}
The pair $(X,Y)$ is log canonical if and only if 
$a_i\geq 0$ for every $i$. If $(X,Y)$ is log canonical
on an open subset containing $W$, then 
$$\mld(W;X,Y)=\min_{\pi(E_i)\subseteq W}\{a_i\}.$$ 
\end{proposition}

\smallskip

The following conjecture is a precise form of
the Inversion of Adjunction Conjecture,
due to Koll\'{a}r and Shokurov. 

\begin{conjecture}
Consider a pair $(X, Y)$ as above and
a normal effective Cartier divisor $D$ on $X$,
such that $D\not\subseteq \bigcup_i Y_i$.
For every nonempty, proper closed subset
$W\subset D$, we have
$$\mld(W;X,D+Y)=\mld(W;D,Y\vert_D),$$
where $Y\vert_D=\sum_iq_i\cdot(Y_i\cap D)$.
\end{conjecture} 

We refer to \cite{kollar} and \cite{kollar2}
for motivation and for a discussion of known results.
Let us mention that there is a more general conjecture
regarding Inversion of Adjunction for
log canonical centers (see, for example, \cite{ambro3}). 
We will prove the following result in Section 3, as an application
of our description of minimal log discrepancies in terms of jet schemes.

\begin{theorem}\label{inv_adj}
The above conjecture is true if $X$ is smooth and $Y=\sum_iq_i\cdot Y_i$,
where $q_i\geq 0$ for all $i$.
\end{theorem}

\begin{remark}
The statement in Theorem~\ref{inv_adj} has been proved in
\cite{ambro1} when $X=\AA^n$, $D$ is a non-degenerate hypersurface,
and $W$ is the origin.
\end{remark}

\begin{corollary}\label{consequence_log_canonical}
Let $(X,Y)$ be a pair as in Theorem~\ref{inv_adj}, with $X$ smooth.
If $D\subset X$ is a normal divisor, such that $D\not\subseteq\bigcup Y_i$,
then
$(X, D+Y)$ is log canonical around $D$
if and only if $(D,Y\vert_D)$
is log canonical.
\end{corollary}

\begin{proof}
Apply Theorem~\ref{inv_adj} with $W=D_{\rm sing}\cup\bigcup_i(Y_i\cap D)$.
Since we have 
$$\mld(D;D,Y\vert_D)=\min\{\mld(W;D,Y\vert_D),1\},$$
$$\mld(D;X,D+Y)=\min\{\mld(W;X,D+Y),0\},$$
we are done since $(X,D+Y)$ is log canonical around $D$
if and only if $\mld(D;X,D+Y)\geq 0$.
\end{proof}

\smallskip

We will consider also the following version of minimal log discrepancy.

\begin{definition}
With the notation in Definition~\ref{def1}, if $W\subset X$
 is a proper irreducible closed subset with generic
point $\eta_W$, then the minimal log discrepency of $(X,Y)$
at $\eta_W$ is
$$\mld(\eta_W;X,Y):=\inf_{c_X(E)=W}a(E;X,Y).$$
If $W=X$, then we put $\mld(\eta_W;X,Y)=0$.
\end{definition}

We collect in the following proposition the basic properties
of this invariant.

\begin{proposition}\label{some_properties}
Let $(X,Y)$ be a pair as above, and
let $W\subset X$ be a proper irreducible closed subset.
\item{\rm (i)} $\mld(\eta_W;X,Y)\geq\mld(W;X,Y)$.
\item{\rm (ii)} If $U\subseteq X$ is open, and $U\cap W\neq\emptyset$,
then $\mld(\eta_W;X,Y)=
\mld(\eta_W;U,Y\vert_U).$
\item{\rm (iii)} $\mld(\eta_W;X,Y)\geq 0$ 
if and only if there is an open subset
$U\subseteq X$ with $U\cap W\neq\emptyset$, and such that $(U,Y\vert_U)$
is log canonical. If $\codim(W,X)\geq 2$, then $\mld(\eta_W;X,Y)\geq 0$
if and only if $\mld(\eta_W;X,Y)\neq -\infty$.
\item{\rm (iv)} If 
$(X,Y)$ is log canonical, and if $\pi : X'\longrightarrow X$ 
is a resolution as in Proposition~\ref{computation}, then
$$\mld(\eta_Z;X,Y)=\inf_{\pi(E_i)=W}\{a_i\}.$$
\item{\rm (v)} There is an open 
subset $U\subseteq X$ such that $U\cap W\neq \emptyset$,
and 
$$\mld(\eta_W;X,Y)=\mld(U\cap W;U,Y\vert_U).$$
\end{proposition}

\section{Log discrepancies and jet schemes}

For the basic definitions and properties of jet schemes, we refer
to \cite{denef} (see also \cite{mustata1} or \cite{mustata2}). 
In particular, we will use freely the construction of motivic
integrals from \cite{denef}.
Recall our context:
$X$ is a normal, $d$-dimensional $\QQ$-Gorenstein variety.
Fix $r\in\NN^*$ such that $r K_X$ is Cartier.

We denote the $m$th jet scheme of $X$ by $X_m$,
and the space of arcs by $X_{\infty}$.  
We have canonical morphisms
$\psi_m : X_{\infty}\longrightarrow X_m$ and $\phi_m : X_m
\longrightarrow X$. When the variety we consider is not obvious, we will
write $\psi_m^X$ and $\phi_m^X$. For every $m$, $j\in\NN$, we put
$X_{m,j}={\rm Im}(X_{m+j}\longrightarrow X_m)$ and similarly
$X_{m,\infty}={\rm Im}(\psi_m)$. It is a theorem of Greenberg from \cite{Gr}
that if $j\gg 0$, then $X_{m,\infty}=X_{m,j}$. In particular,
$X_{m,\infty}$ is constructible.

On each jet scheme $X_m$ there is an $\AA^1$-action $\bullet$ such that
the natural projections are compatible with these actions. Here $\AA^1$
is considered as a monoidal scheme under usual multiplication.
Moreover, there are ``zero-sections'' $\sigma_m : X\longrightarrow X_m$ such that
for every $\gamma\in X_m$, we have $0\bullet\gamma=\sigma_m(\phi_m(\gamma))$.
Note that if $T\subseteq X_m$ is invariant
under the $\AA^1$-action, then $\phi_m(T)=\sigma_m^{-1}(T)$, hence it is closed
in $X$. For more details we refer to \cite{mustata2}.

We introduce now two subschemes of $X$ which measure its singularities.
Let $i : U=X_{\rm reg}\longrightarrow X$
be the open immersion corresponding to the smooth part of $X$.
We have a canonical morphism
$$(\Omega_X^d)^{\otimes r}\longrightarrow i_*(\Omega_U^d)^{\otimes r}=
\cO_X(rK_X).$$ 
This defines a closed subscheme $Z\subset X$ of ideal $\cI_Z$, such that
the image of the above morphism is $\cI_Z\otimes\cO_X(rK_X)$. 

We consider also the Jacobian subscheme $Z'$ defined by the 
Jacobian ideal $\cI_{Z'}:={\rm Fitt}_d(\Omega_X^1)$ (the $d$th Fitting ideal
of $\Omega_X^1$). Working locally, we may assume that $X\subset \AA^N$
is defined by $(f_i)_i$. Then $\cI_{Z'}$ is generated by the restrictions
to $X$ of the $(N-d)$ minors of the Jacobian matrix 
$(\partial f_i/\partial X_j)_{i,j}$.

It is clear that we have ${\rm Supp}(Z)\subseteq X_{\rm sing}={\rm Supp}(Z')$.
Moreover, if $X$ is locally complete intersection, then 
we may take $r=1$ and in this case $\cI_Z=\cI_{Z'}$.

Recall that to every closed subscheme $T\hookrightarrow X$
we have an associated 
function $F_T : X_{\infty}\longrightarrow \NN\cup\{\infty\}$
which measures the order of vanishing of an arc along $T$.
For every $e\in\NN$, let $X^{(e)}_{\infty}=F_{Z'}^{-1}(e)\subseteq X_{\infty}$.
Similarly, if $m\geq e$, we denote by $X_m^{(e)}$ the
set of jets in $X_m$ vanishing along $I_{Z'}$ with order exactly $e$.
We put also $X^{(e)}_{m,j}=X_m^{(e)}\cap X_{m,j}$, and similarly 
for $X_{m,\infty}^{(e)}$.

It follows from Lemma~4.1 in \cite{denef} 
and its proof (see also \cite{looijenga}) that if $m\geq e$,
then $X_{m+1,\infty}^{(e)}\longrightarrow X_{m,\infty}^{(e)}$
is locally trivial with fiber $\AA^d$. Moreover, $X_{m,\infty}^{(e)}$
is a locally closed subset of $X^{(e)}_m$.

Consider for all closed subschemes $T\subseteq X$ the corrresponding
subsets $\psi_m^{-1}(T_m)$, and let $A$ be an element in the algebra
generated by all these subsets. We define $\codim(A)$ as follows. 
Suppose first that $A\subseteq X_{\infty}^{(e)}$ for some $e$,
and let us write $A=\psi_m^{-1}(B)$, where we may take $m\geq e$.
We put $\codim(A):=(m+1)d-\dim(B\cap X_{m,\infty}^{(e)})$,
and since $\dim\psi_m(X_{\infty})=(m+1)d$ (see Lemma~4.3 in \cite{denef}),
it follows that $\codim(A)$ is a nonnegative integer.
Moreover, the result mentioned in the above paragraph
shows that the definition does not depend on which $m$ we have chosen.
In general, we put $\codim(A):=\min_{e\in\NN}\codim(A\cap X_{\infty}^{(e)})$,
with the convention that $\codim(\emptyset)=\infty$. Note that 
$\codim(A)=\infty$ if and only if $A\subseteq Z'_{\infty}$.

Let $\mu(A)$ be the Hodge realization of the motivic
measure of $A$ (see \cite{denef}). $\mu(A)$ is a Laurent
power series in two variables $u^{-1}$ and $v^{-1}$.
If $\mu(A)\neq 0$, then there is only one monomial in $\mu(A)$
of maximal degree, namely
 $c(uv)^{-\codim(A)}$, where $c$ is a positive integer.
Moreover, $\mu(A)=0$ if and only if $\codim(A)=\infty$.

\begin{remark}
Suppose that $A$ is a finite intersection of sets of the form 
$F_{T_i}^{-1}(\geq m_i)$. One can show that in this case $\codim(A)<\infty$
and, in fact, $\codim(A)=\codim(\psi_m(A),\psi_m(X_{\infty}))$, if $m\gg 0$.
If $X$ is smooth, then one can also check that $\codim(A)$
is the codimension of $A$ as a closed subset of $X_{\infty}$,
but we do not know if this remains true in general.
\end{remark} 

\bigskip

We give now the version of the Change of Variable formula we will need
(this is essentially the same version that was used in \cite{yasuda}).
Let $X$ be as above, $Z$ the subscheme we have defined, and
$W\subseteq X$ a closed subset. We consider $Y=\sum_{i=1}^k q_i\cdot Y_i$,
where $q_i\in\RR$, and $Y_i\subset X$ are proper closed subschemes. 
For $e\in\NN$, $m=(m_i)_i\in\NN^k$, we put $A=\bigcap_i F_{Y_i}^{-1}(m_i)\cap 
F_Z^{-1}(e)\cap\psi_0^{-1}(W)$.

\begin{theorem}\label{change_of_var}
If $\pi : X'\longrightarrow X$ is a proper, birational morphism, with
$X'$ smooth, then
$$\int_A(uv)^{(1/r)F_Z}=\int_{\pi_{\infty}^{-1}(A)}
(uv)^{-F_{K_{X'/X}}}.$$
\end{theorem}

\begin{proof}
By the Change of Variable formula in \cite{denef}, we have
$$\int_A(uv)^{(1/r)F_Z}=\int_{\pi_{\infty}^{-1}(A)}
(uv)^{(1/r)F_Z\circ\pi_{\infty}-F_{Z''}},$$
where $\cI_{Z''}$ is the sheaf of ideals such that 
$\pi^*\Omega_X^d\longrightarrow \cI_{Z''}\otimes\Omega_{X'}^d$
is an epimorphism. It follows from definition that
$$\pi^{-1}(\cI_Z)\cdot\cO(-rK_{X'/X})=\cI_{Z''}^r.$$
Since 
$F_Z\circ\pi_{\infty}=F_{\pi^{-1}(Z)}$, we deduce the formula in the statement.
\end{proof}

The following is the main technical ingredient, which will allow
us to connect log discrepancies and jet schemes. We first fix the
notation. Let $(X,Y)$ and $Z\subset X$ be as above. We consider also
a closed nonempty
subset $W\subset X$. Fix a proper, birational morphism
$\pi : X'\longrightarrow X$, such that $X'$ is smooth, and such that
$\pi^{-1}(Y)\cup\pi^{-1}(Z)\cup{\rm Ex}(\pi)$
is a divisor with simple normal crossings. 
Recall that $\pi^{-1}(Y):=
\sum_iq_i\pi^{-1}(Y_i)$. If $W\neq X$, then we put also the
condition that $\pi^{-1}(W)$, together with the above union, is a divisor
with simple normal crossings. We write $\pi^{-1}(Y_i)=\sum_{j=1}^sy_{i,j}D_j$,
$\pi^{-1}(Z)=\sum_{j=1}^sz_jD_j$, and $K_{X'/X}=\sum_{j=1}^sk_jD_j$. 

\begin{theorem}\label{formula_for_dim}
With the above notation, for every $e\in\NN$, $m\in\NN^k$, we have
$$\codim\left(\bigcap_iF_{Y_i}^{-1}(m_i)\cap 
F_Z^{-1}(e)\cap\psi_0^{-1}(W)\right)=
\frac{e}{r}+\min_{\nu}\sum_{j=1}^s(k_j+1)\nu_j.$$
Here the infimum is over all $\nu\in\NN^s$ with $\bigcap_{\nu_j\neq 0}
D_j\neq\emptyset$, and such that
$\sum_{j=1}^s\nu_jy_{i,j}=m_i$ for all $i$, and $\sum_{j=1}^s\nu_jz_j=e$.
If $W\neq X$, then we have to add also the condition that there is
$\nu_j\neq 0$ such that $\pi(D_j)\subseteq W$.
By convention, the minimum over an empty set is $\infty$.
\end{theorem}

\begin{proof}
Let $A=\bigcap_iF_{Y_i}^{-1}(m_i)\cap 
F_Z^{-1}(e)\cap\psi_0^{-1}(W)$. By definition,
$$\int_A(uv)^{(1/r)F_Z}=\mu(A)(uv)^{e/r}.$$
We treat only the case $\mu(A)\neq 0$, the changes for the other case
being obvious. Therefore the integral is given by a
Laurent power series in $u^{-1}$ and $v^{-1}$ (with rational exponents)
of degree $2((e/r)-\codim(A))$.

On the other hand, a direct computation shows that
$$\int_{\pi_{\infty}^{-1}(A)}(uv)^{-F_{K_{X'/X}}}
=\sum_{\nu}E(D_{\nu}^{\circ})(uv-1)^{|\nu|}
(uv)^{-d-\sum_j(k_j+1)\nu_j},$$
where the sum is over all $\nu\in\NN^s$ such that $\sum_i\nu_jy_{i,j}=m_i$,
for all $i$,
and $\sum_j\nu_jz_j=e$. We have put
$|\nu|={\rm Card}\{j\vert\nu_j>0\}$. 
If $W\neq X$, then we have to add also the condition
that there is at least one $\nu_j>0$ such that
$\pi(D_j)\subseteq W$. We have used the notation $D^{\circ}_{\nu}
=\bigcap_{\nu_j\neq 0}
D_j\setminus\bigcup_{\nu_j=0}D_j$. The proof of this formula follows along the
same lines as the proof of Theorem~2.15 in \cite{craw}.

We deduce that the degree of the integral over $\pi_{\infty}^{-1}(A)$
is equal to 
$-2(\min_{\nu}\sum_j(k_j+1)\nu_j)$, where $\nu$ runs over the set in the statement
of the theorem. By Theorem~\ref{change_of_var} the two integrals are equal,
and comparing their degrees we get the formula for $\codim(A)$.
\end{proof}

\begin{remark}
When $X$ is smooth, $W=X$ and $Y$ is a hypersurface, the formula in
Theorem~\ref{formula_for_dim} follows also from the computation of
motivic Igusa zeta function in \cite{denef2}. Under the same
assumption on $X$ and $W$, but for $Y=q_1\cdot Y_1-q_2\cdot Y_2$,
this is contained in \cite{elm}.
\end{remark}

\begin{corollary}\label{formula_for_dimension10}
With the notation in Theorem~\ref{formula_for_dim}, we have
$$\codim\left(\bigcap_i F_{Y_i}^{-1}(\geq m_i)
\cap F_Z^{-1}(\geq e)\cap\psi_0^{-1}(W)\right)=
\min_{\nu}\sum_{j=1}^s\left(\frac {z_j}{r}+k_j+1\right)\nu_j.$$
Here the infimum is over all $\nu\in\NN^s$ with $\bigcap_{\nu_j\neq 0}
D_j\neq\emptyset$, and such that
$\sum_{j=1}^s\nu_jy_{i,j}\geq m_i$
for all $i$, and $\sum_{j=1}^s\nu_jz_j\geq e$.
If $W\neq X$, then we have to add also the condition that there is
$\nu_j\neq 0$ such that $\pi(D_j)\subseteq W$.
\end{corollary}

\begin{proof}
If $m\in\NN^k$ and $e\in\NN$, then
we put $A_{m,e}$ for the set in Theorem~\ref{formula_for_dim}, and 
$A'_{m,e}$ for the set in the above statement. We put $m'\geq m$
if $m'_i\geq m_i$ for all $i$.
We have
$$\bigcup_{m'\geq m,e'\geq e}A_{m',e'}\subseteq A'_{m,e},$$
and the complement lies inside $\bigcup_i(Y_i)_{\infty}\cup Z_{\infty}$.
This gives 
$$\codim(A'_{m,e})=\min_{m',e'}\codim(A_{m',e'}),$$
and we conclude by Theorem~\ref{formula_for_dim}.
\end{proof}

We give now our characterization of minimal log discrepancies in terms
of spaces of arcs. Using Proposition~\ref{basic_properties}(ii), 
it is easy to see that for the computation of $\mld(W;X,Y)$ we 
may assume that $W\neq X$.

\begin{theorem}\label{char_log_canonical}
Let $(X,Y)$ be a pair as above, with $Y=\sum_{i=1}^kq_i\cdot Y_i$.
If $W\subset X$ is a proper closed subset, and if $\tau\in{\mathbb R}_+$,
then the following are equivalent:
\item{\rm (i)} $\mld(W;X,Y)\geq\tau$.
\item{\rm (ii)} For every $e\in\NN$, $m\in\NN^k$, we have
\begin{equation}\label{condition_on_jets2}
\codim\left(\bigcap_{i=1}^kF_{Y_i}^{-1}(m_i)\cap F_Z^{-1}(e)
\cap\psi_0^{-1}(W)\right)\geq\frac{e}{r}+\sum_{i=1}^k
q_im_i+\tau.
\end{equation}

If $q_i\geq 0$ for all $i$, then the above conditions are also
equivalent with
\item{\rm (iii)} For every $e\in\NN$, $m\in\NN^k$, we have
\begin{equation}\label{formula3}
\codim\left(\bigcap_{i=1}^k F_{Y_i}^{-1}(\geq m_i)\cap F_Z^{-1}(\geq e)
\cap\psi_0^{-1}(W)\right)
\geq\frac{e}{r}+\sum_{i=1}^kq_i m_i+\tau.
\end{equation} 

Moreover, if $\pi$ is a resolution 
as in Theorem~\ref{formula_for_dim},
then in {\rm (ii)} and {\rm (iii)} above it is
enough to put the conditions in {\rm (\ref{condition_on_jets2})} and 
{\rm (\ref{formula3})},
respectively, only for finitely many $e$ and $m$, depending on the
numerical data of the resolution.
\end{theorem}

\begin{proof}
Let $\pi$ be a resolution as in Theorem~\ref{formula_for_dim}. We keep the
notation in that theorem. By restricting to a suitable open neighbourhood
of $W$, we may assume that $\pi(D_j)\cap W\neq\emptyset$, for all $j$.
In this case Proposition~\ref{computation} shows that
$\mld(W;X,Y)\geq\tau$ if and only if $k_j+1-\sum_iq_iy_{i,j}\geq 0$
for all $j$, and
$k_j+1-\sum_iq_iy_{i,j}\geq\tau$ for all $j$ such that $\pi(D_j)\subseteq W$.

Suppose first that $\mld(W;X,Y)\geq\tau$.
If $\nu\in\NN^s$ is such that
$\sum_j\nu_jy_{i,j}=m_i$, and such that $\nu_l\geq 1$
for some $l$ with $\pi(D_l)\subseteq W$, 
then $\sum_j(k_j+1)\nu_j\geq\sum_iq_im_i+\tau$, and we deduce
(\ref{condition_on_jets2}) from the formula in Theorem~\ref{formula_for_dim}.
This proves (i)$\Rightarrow$(ii), and 
(i)$\Rightarrow$(iii) follows similarly, using 
Corollary~\ref{formula_for_dimension10}.

We show now (ii)$\Rightarrow$(i). We take first $j$
such that $\pi(D_j)\subseteq W$.
If $k_j+1<\sum_iq_iy_{i,j}+\tau$, take $m\in\NN^k$ given by 
$m_i=y_{i,j}$ for all $i$, and let $e=z_j$. 
By taking $\nu_j=1$,
and $\nu_{j'}=0$ for $j'\neq j$, the formula in Theorem~\ref{formula_for_dim}
gives
$$\codim\left(\bigcap_{i=1}^k F_{Y_i}^{-1}(m_i)
\cap F_Z^{-1}(e)\cap \psi_0^{-1}(W)\right)
<\frac{e}{r}+\sum_{i=1}^kq_im_i+\tau,$$
a contradiction with (ii). 

We take now $j$ such that 
$\pi(D_j)\not\subseteq W$. Since 
$\pi(D_j)\cap W\neq\emptyset$, 
there is $j'$ such that $\pi(D_{j'})\subseteq W$,
and $D_j\cap D_{j'}\neq\emptyset$. For every $\alpha\in\NN$, we take 
$\nu\in\NN^s$ such that $\nu_j=\alpha$, $\nu_{j'}=1$, and $\nu_{j''}=0$
if $j''\neq j$, $j'$. If $m\in\NN^k$ is such that $m_i=y_{i,j}\alpha+
y_{i,j'}$, and if $e=z_j\alpha+z_{j'}$, then it follows from 
(\ref{condition_on_jets2}) 
and the formula in Theorem~\ref{formula_for_dim} that
$$(k_j+1)\alpha+(k_{j'}+1)\geq\sum_{i=1}^kq_i(y_{i,j}\alpha+y_{i,j'})+\tau.$$
It is clear that there is a value for $\alpha$, depending on the numerical
data of the resolution, such that the above inequality implies
$k_j+1\geq\sum_iq_iy_{i,j}$. We have thus shown
that $\mld(W;X,Y)\geq\tau$.

Note that (iii) trivially implies (ii), 
as the codimension in (\ref{condition_on_jets2})
is always greater or equal to the codimension in (\ref{formula3}).  
As the last assertion in the theorem follows from the above arguments,
we are done.
\end{proof}

\begin{remark}\label{char_log_canonical2}
We can deduce from the above theorem a condition 
for $(X,Y)$ to be log canonical, in terms of arcs. Namely, $(X,Y)$
is log canonical if and only if for every $e\in\NN$ and every 
$m\in\NN^k$, we have
\begin{equation}
\codim\left(\bigcap_{i=1}^kF_{Y_i}^{-1}\cap F_Z^{-1}(e)\right)
\geq\frac{e}{r}+\sum_{i=1}^kq_im_i.
\end{equation}
Moreover, if $q_i\geq 0$ for all $i$, then the above condition is equivalent
with
\begin{equation}
\codim\left(\bigcap_iF_{Y_i}^{-1}(\geq m_i)\cap F_Z^{-1}(\geq e)\right)
\geq\frac{e}{r}+\sum_iq_im_i,
\end{equation}
for every $e\in\NN$ and every $m\in\NN^k$. In order to see this,
it is enough to apply the above theorem for $W=\bigcup_i{\rm Supp}(Y_i)\cup
X_{\rm sing}$.

Note that this characterization of log canonical singularities was proved
in the case when ${\rm Supp}(Z)\subseteq\bigcup_i{\rm Supp}(Y_i)$ in 
\cite{yasuda}. The above proof of Theorem~\ref{char_log_canonical}
is inspired from his proof.
\end{remark}

\begin{remark}\label{classical}
One can give an analogous description of minimal log discrepancies
in the usual setting of Mori Theory. Suppose that $X$ is a 
$d$-dimensional normal variety,
and that $D$ is a $\QQ$-divisor on $X$ such that $r(K_X+D)$ is Cartier
for some positive integer $r$. For simplicity, we assume that 
$D$ is effective, so we have a canonical morphism
$$(\wedge^d\Omega_X)^{\otimes r}
\longrightarrow(\wedge^d\Omega_X)^{\otimes r}\otimes\cO_X(rD)
\longrightarrow\cO_X(r(K_X+D)).$$
We have a closed subscheme $T\subseteq X$ defined by the ideal ${\mathcal I}_T$,
such that the image of the above composition is 
${\mathcal I}_T\otimes\cO_X(r(K_X+D))$. 
Note that in this case, this scheme depends also on $D$.

The same arguments as above show, for example, that if $\tau\in\RR_+$, then
$\mld(W;X,D)\geq\tau$ if and only if 
$$\codim(F_T^{-1}(e)\cap\psi_0^{-1}(W))\geq \frac{e}{r}+\tau$$
for every $e\in\NN$. 
\end{remark}

\section{Inversion of Adjunction}

In the case of a hypersurface, it is easy to understand
the set of jets which can be lifted to the arc space.
This will be enough to give a proof of 
Theorem~\ref{inv_adj}. 

Let us fix the notation for this section.
We consider a smooth variety $X$, with $\dim\,X=d$, and a
divisor $D\subset X$ which is irreducible and reduced.
Let $Z\subset D$ be the jacobian subscheme of $D$ defined by the ideal
$I_Z$. Recall that $D_{m,e}={\rm Im}(D_{m+e}\longrightarrow D_m)$,
while $D_{m,\infty}={\rm Im}(D_{\infty}\longrightarrow D_m)$. 
In addition, if we restrict to those jets with order $e$ along $I_Z$,
then we put $(e)$ as a superscript.

\begin{lemma}\label{equations_codim1}
Given $D$ as above, and $m$, $e\in\NN$, with $m\geq e$,
we have $D^{(e)}_{m,\infty}=D^{(e)}_{m,e}$.
Moreover, if $\eta : X_{m+e}\longrightarrow X_m$ is the 
canonical projection, then
$\eta^{-1}(D_{m,e}^{(e)})=D_{m+e}^{(e)}$. 
\end{lemma}

\begin{proof}
We have to show that if $u\in D^{(e)}_{m+e}$, then there is 
$v\in D_{\infty}$ such that $u$ and $v$ have the same image in $D_m$.
In addition, if $w\in X_{m+e}$ is such that $\eta(u)=\eta(w)$, then $w\in D_{m+e}$.

Let $u_0=\phi^D_m(u)$. By restricting to an open neighbourhood of $u_0$,
we may assume that we have a regular system of parameters at $u_0$,
denoted by $x_1,\ldots,x_d$. We may also assume that $D$ is defined by 
an equation $f$. Note that the regular system of parameters induces
an isomorphism $\widehat{\cO}_{X,u_0}\simeq\CC[[T_1,\ldots,T_d]]$, and we will
identify $f$ with a power series via this isomorphism. 

For every $p$, we have an isomorphism
$$(\phi_p^X)^{-1}(u_0)\simeq (t\CC[t]/(t^{p+1}))^d,$$
which maps a morphism $\gamma$ to $(\gamma_i)_i$,
where $\gamma_i=\gamma(x_i)$. Note that $\gamma\in D_p$
if and only if $f(\gamma_1,\ldots,\gamma_d)=0$. Similar
considerations apply when $p=\infty$.

In order to finish the proof, it is enough to prove the following assertions.
Suppose that $\gamma\in (t\CC[[t]])^n$ is such that $\ord\,f(\gamma)\geq m+e+1$
and 
$\ord(\partial f/\partial T_1(\gamma),\ldots,\partial f/\partial T_n(\gamma))=e$.
Then there is
$\delta\in (\CC[[t]])^n$, such that $f(\gamma+t^{m+1}\delta)=0$.
Moreover, if $\gamma'\in (\CC[[t]])^n$, then $\ord\,f(\gamma+t^{m+1}\gamma')
\geq m+e+1$.

Consider the  Taylor expansion:
$$f(\gamma+t^{m+1}\gamma')=f(\gamma)+t^{m+1}\sum_{i=1}^d
\frac{\partial f}{\partial T_i}(\gamma)
\cdot\gamma'_i+t^{2m+2}\cdot(\ldots).$$
Since $\ord\,f(\gamma)\geq m+e+1$, $\ord\frac{\partial f}{\partial T_i}(\gamma)
\geq e$, and $m\geq e$, we deduce $\ord\,f(\gamma+t^{m+1}\gamma')\geq
m+e+1$.
This gives the second of the above assertions. Moreover, it is easy to
see
from the above formula that there is $\delta$ such
that
$f(\gamma+t^{m+1}\delta)=0$. In fact, the terms of order zero in
$\delta$
are obtained solving a linear equation, while the higher terms can be
deduced
by a recursive argument. Hence we get the first of the above assertions.
Alternatively, this statement 
can be deduced also from Newton's Lemma (see \cite{Gr}).
\end{proof}

Consider now the following situation. Let $m\in\NN$, and
let also $R\subseteq X_m$ be an irreducible closed subset which is
invariant under the $\AA^1$-action on $X_m$. 
Suppose that $a\leq m$ is such that ${\rm ord}\,\gamma(I_D)\geq a$
for every $\gamma\in R$, where $I_D$ is the ideal defining $D$ in $X$.
We assume that $\phi_m^X(R)\cap D\neq\emptyset$, and
we put $S=R\cap D_{m,\infty}$ and ${\mathcal S}=(\psi_m^D)^{-1}(S)$.
Let
$$e=\min\{\ord\,\gamma(I_Z)\vert\gamma\in {\mathcal S}\}.$$

\begin{lemma}\label{prelim}
With the above notation, we have
\begin{enumerate}
\item $S$ is non-empty and $e<\infty$.
\item If ${\mathcal S}^{\circ}
:=\{\gamma\in {\mathcal S}\vert\ord\,\gamma(I_Z)=e\}$ (which is open in 
${\mathcal S}$), 
then 
$$\codim({\mathcal S}^{\circ},D_{\infty})\leq\codim(R,X_m)+e-a.$$
\end{enumerate}
\end{lemma} 

\begin{proof}
Let $x\in\phi_m^X(R)\cap D$. We denote by $x_m$ the image of $x$ by
the zero section to $D_m$. Using the $\AA^1$-action on $D_m$,
we see that $x_m\in S$. 
The second assertion in (1) can be proved as follows.
Let $\mu : D'\longrightarrow D$ be a resolution of
singularities for $D$. The set $f_{\infty}^{-1}({\mathcal S})$
is nonempty as it contains the zero section over any point in $f^{-1}(x)$.
Since it is the inverse image of a closed subset in $D'_m$, 
and since $D'$ is smooth, it can not be
contained in $f^{-1}(Z)_{\infty}$ (see, for example, 
Corollary 3.8 in \cite{mustata2}). This shows that ${\mathcal S}\not\subseteq
Z_{\infty}$.

Fix now $p\geq\max\{m,e\}$. Let ${\mathcal R}:=(\psi_m^X)^{-1}(R)$
We denote by ${\mathcal R}_{p+e}$
and ${\mathcal R}_p$ the projections of ${\mathcal R}$ to
$X_{p+e}$ and $X_p$, respectively. 
Similarly, let ${\mathcal S}^{\circ}_p$
be the projection of ${\mathcal S}^{\circ}$ to $X_p$. 
We denote by $g :{\mathcal R}_{p+e}\longrightarrow{\mathcal R}_p$
the canonical projection.

Let $T={\mathcal R}_{p+e}\cap D_{p+e}$. If $\gamma\in T$ has order $e'\leq e$
along $I_Z$, then Lemma~\ref{equations_codim1} 
shows that $g(\gamma)$ lies over $S$. Hence $e'=e$.
Therefore the set $T^{\circ}
:=\{\gamma\in T\vert{\rm ord}\gamma(I_Z)=e\}$ is an open subset of
$T$. Again, Lemma~\ref{equations_codim1} implies that $g$ induces a surjective map
$f : T^{\circ}\longrightarrow {\mathcal S}^{\circ}_p$. 

It is clear that all the fibers of $f$ have dimension at most $de$.
On the other hand, $T$ is cut out in ${\mathcal R}_{p+e}$ by
$p+e-a+1$ equations. If we put $r=\codim(R,X_m)$, then 
every irreducible component of $T$ has dimension at least
$(p+1)(d-1)+de-(r+e-a)$. Therefore $\dim{\mathcal S}_p^{\circ}
\geq (p+1)(d-1)-(r+e-a)$, hence $\codim({\mathcal S}^{\circ},D_{\infty})
\leq r+e-a$.
\end{proof}

\smallskip

We can prove now the case of Inversion of Adjunction which was 
stated in Section 1.

\begin{proof}[Proof of Theorem~\ref{inv_adj}]
The inequality 
$$\mld(W;X,D+Y)\leq\mld(W;D,Y\vert_D)$$
is well-known in general and follows by adjunction (see, for example,
the proof of Proposition~7.3.2 in \cite{kollar}). We recall the 
argument for completeness.

Let $\pi : X'\longrightarrow X$ be proper, birational, such that $X'$
is smooth, and $\pi^{-1}(Y)\cup\pi^{-1}(D)\cup\pi^{-1}(W)$ is a divisor
with simple normal crossings. Write 
$$K_{X'/X}-\pi^{-1}(Y)=\sum_i(a_i-1)E_i,$$ 
and $\pi^{-1}(D)=\widetilde{D}+\sum_ib_iE_i$.
Note that by hypothesis, the strict transform $\widetilde{D}$
of $D$ does not appear in $K_{X'/X}-\pi^{-1}(Y)$.
 
The restriction $\pi_0 : \widetilde{D}\longrightarrow D$
is a log resolution of $(D,Y\vert_D\cup W)$, and the adjunction formula gives
$$K_{\widetilde{D}/D}-\pi_0^{-1}(Y\vert_D)=\sum_i(a_i-b_i-1)
E_i\vert_{\widetilde{D}}.$$
Since $\pi(\widetilde{D})\not\subseteq W$, we see that 
if $\widetilde{D}\cap E_i\neq\emptyset$, then $\pi_0(\widetilde{D}\cap E_i)
\subseteq W$ if and only if $\pi(E_i)\subseteq W$. This is enough to give the 
inequality we have claimed. The reverse inequality is not obvious, 
as some of the divisors $E_i$ might not intersect $\widetilde{D}$.

We turn now to the proof of this reverse inequality.
Suppose that 
$\mld(W;X,D+Y)<\tau$, for some $\tau\in\RR_+$.
It follows from Theorem~\ref{char_log_canonical}
applied to the smooth variety $X$ 
that there are $m$, $a\in\NN$, $b\in\NN^k$,  such that $m\geq\max\{a,b_i\}$, and if
$$A=\{\gamma\in X_m\vert
\ord\,\gamma(I_D)\geq a,\ord\,\gamma(I_{Y_i})\geq b_i,\ord\gamma(I_W)\geq 1\},$$
then there is an irreducible component $R$ of $A$, with 
$\codim(R,X_m)<a+\sum_ib_iq_i+\tau$.

It is clear that
$R$ satisfies the hypothesis of Lemma~\ref{prelim}.
With notation as in the lemma, we get the subset
$${\mathcal S}^{\circ}\subseteq D_{\infty}\cap F_Z^{-1}(e)\cap\bigcap_i
F_{Y_i}^{-1}(\geq b_i),$$
such that $\codim({\mathcal S}^{\circ},D_{\infty})<\sum_iq_ib_i+e+\tau$.
Theorem~\ref{char_log_canonical} shows that
$\mld(W;D,Y\vert_D)<\tau$.
This completes the proof of the theorem.
\end{proof}

We apply now Theorem~\ref{inv_adj} and~\ref{char_log_canonical} to deduce a 
characterization of terminal hypersurfaces. Fix a divisor 
$D\subset X$, where $X$ is smooth of dimension $d$,
and $D$ is reduced and irreducible.
Recall that by Theorem~3.3
in \cite{mustata2}, $D$ has canonical (or equivalently, rational)
singularities if and only if $D_m$ is irreducible for every $m$. Moreover, it is 
shown in \cite{mustata2} that in this case $D_m$ is a locally complete intersection
variety, of dimension $(m+1)(d-1)$. The following result similarly
characterizes terminal singularities, giving a positive answer to a question
of Mirel Caib\v{a}r.

\begin{theorem}\label{terminal}
If $D\subset X$ is an irreducible and reduced divisor on a smooth
variety $X$, then $D$
has terminal singularities if and only if $D_m$ is normal for every $m\in\NN$.
\end{theorem}

\begin{proof}
By the results in \cite{mustata2}, 
we may assume that $D_m$ is a locally complete intersection
variety of dimension $(m+1)(d-1)$. Therefore $D_m$ is normal if and only if
$\dim(D_m)_{\rm sing}\leq (m+1)(d-1)-2$. 

Moreover, if $\phi_m :D_m\longrightarrow D$ is the canonical projection, then
$(D_m)_{\rm sing}=\phi_m^{-1}(D_{\rm sing})$. Indeed, the inclusion
``$\subseteq$'' is trivial. To see the reverse inclusion, note that
$D_m\subset X_m$ is defined by $(m+1)$ equations. Since
$\dim\,D_m=\dim\,X_m-(m+1)$, if $u\in D_m$ is a smooth point, then
these $(m+1)$ equations are part of a regular system of parameters at $u$.
As the scheme defined by the first equation is
locally isomorphic to $D\times\AA^{md}$,
$D$ has to be smooth at $\phi_m(u)$.

Therefore $D_m$ is normal for every $m$ if and only if
$\dim\,\phi_m^{-1}(D_{\rm sing})\leq (m+1)(d-1)-2$ for every $m$. 
By Theorem~\ref{char_log_canonical}, 
this is equivalent with $\mld(D_{\rm sing};X,D)\geq 2$.
Since this minimal log discrepancy is an integer, this is further
equivalent to $\mld(D_{\rm sing};X,D)>1$. By Theorem~\ref{inv_adj}, we have
$$\mld(D_{\rm sing};X,D)=\mld(D_{\rm sing};D).$$
 As by definition 
$\mld(D_{\rm sing};D)>1$ if and only if $D$ has terminal singularities, 
we are done.
\end{proof}

\begin{remark}
In fact, the above argument can be used to show
that if $D$ is a normal divisor on a smooth $d$-dimensional variety $X$,
then $D$ has log canonical singularities 
if and only if $D_m$ has pure dimension 
for every $m$. Indeed, $D_m$ has pure dimension if and only if 
$\dim\,D_m=(m+1)\dim\,D$. By Remark~\ref{char_log_canonical2}, this is true
for every $m$ if and only if $(X,D)$ is log canonical. This is equivalent
with $D$ being log canonical by
Corollary~\ref{consequence_log_canonical}.

Suppose now that $D$ is a normal divisor with log canonical
singularities.
The argument in the proof of Theorem~\ref{terminal} shows that $(D_m)_{\rm sing}
=\phi_m^{-1}(D_{\rm sing})$. Moreover, we see that 
$\codim((D_m)_{\rm sing},D_m)\geq \mld(D_{\rm sing};D)$, for all $m$, and
equality is achieved for some $m$. In particular, this implies
Theorem~3.3
in \cite{mustata2}: $D$ has canonical singularities if and only if 
$D_m$ is irreducible for every $m$. 
\end{remark}

\section{Semicontinuity of minimal log discrepancies}

In this section we prove a semicontinuity statement for
minimal log discrepancies in the case of a smooth ambient variety.
Recall the following conjecture from \cite{ambro}.

\begin{conjecture}\label{semicontinuity}
If $X$ is a normal, $\QQ$-Gorenstein variety, 
and if $Y=\sum_iq_i\cdot Y_i$, where $q_i\in\RR_+$ and $Y_i\subset X$ is a proper
closed subscheme, for all $i$, then
the function
$x\in X\longrightarrow\mld(x;X,Y)$ is lower semicontinuous.
\end{conjecture}

It was shown in \cite{ambro}
that this conjecture is equivalent with the following one.

\begin{conjecture}\label{inequality}
Let $X$ and $Y$  be as in Conjecture~\ref{semicontinuity}.
For every two irreducible closed subsets $V\subset W\subset X$, we have
$$\mld(\eta_V;X,Y)\leq\mld(\eta_W;X,Y)+\codim(V,W).$$
\end{conjecture}

\begin{remark}
In fact, in \cite{ambro}, $Y$ is assumed to be a divisor. On the other hand,
all the arguments can be extended to the case of an arbitrary subscheme.

One reason for conjecturing the above statements in \cite{ambro}
was to explain a conjecture of V.~Shokurov from \cite{shokurov},
which was the particular case $W=X$ in Conjecture~\ref{inequality}.
\end{remark}

We will show that Conjecture~\ref{inequality} is true if the ambient variety
is smooth.

\begin{theorem}
Let $X$ be a smooth variety, and $Y=\sum_{i=1}^kq_i\cdot Y_i$, where
$q_i\in\RR_+$ and $Y_i\subset X$ is a proper closed subscheme, for all $i$.
For every two irreducible closed subsets $V\subset W$, we have
$$\mld(\eta_V;X,Y)\leq\mld(\eta_W;X,Y)+\codim(V,W).$$
\end{theorem}

\begin{proof}
By taking a sequence of intermediate subvarieties, it is enough to
consider the case when $\codim(V,W)=1$. If $W=X$, then $V$ is a divisor.
We clearly have $\mld(\eta_V;X,Y)\leq\mld(\eta_V;X)=1$, which
completes this case, as $\mld(\eta_X;X,q\cdot Y)=0$.

From now on, we suppose that $W\neq X$, so $\codim(V,X)\geq 2$.
If $(X,Y)$ is not log canonical on any open subset meeting $V$,
then $\mld(\eta_V;X,Y)=-\infty$, and there is nothing to prove. 
If this is not the case, by restricting to a suitable open subset, we may assume
that $(X,Y)$ is log canonical. Moreover, we may restrict to a suitable
open subset meeting $V$ in order to have
$\mld(V;X,q\cdot Y)=\mld(\eta_V;X,q\cdot Y)$ 
(see Proposition~\ref{some_properties}(v)).
Up to this point, the argument
holds for arbitrary $X$. 

Let $\tau=\mld(\eta_W;X,Y)$. Since $X$ is smooth, by 
Theorem~\ref{char_log_canonical} there exists
$m=(m_i)_i\in\NN^k$ such that
$\codim(A,X_{\infty})\leq \sum_iq_im_i+\tau$, where
$A=\bigcap_iF_{Y_i}^{-1}(\geq m_i)\cap \psi_0^{-1}(W)$.
Fix $p\gg 0$ (depending on $m$), and let $B=\psi_p(A)$, so that
$\dim\,B\geq (p+1)\dim\,X-\sum_iq_im_i-\tau$. 

We claim that we can choose $m$ such that there is an irreducible
component $T$ of $B$ with $\dim\,T\geq (p+1)\dim\,X-\sum_iq_im_i-\tau$,
and such that $\phi_p(T)=W$. Indeed, note first that for every irreducible
component $T$, $\phi_p(T)$ is a closed subset of $W$. This follows since
$T$ is invariant under
the $\AA^1$-action on $X_m$. If we can not
find $T$ as claimed, then we may restrict to a suitable open subset meeting $W$
to deduce $\mld(\eta_W;X,Y)>\tau$, a contradiction.
 For this we use the fact that by Theorem~\ref{char_log_canonical},
in order to compute minimal log discrepancies it is enough to check finitely
many jet schemes, depending on a log resolution of $(X,Y\cup W)$.
Therefore we can find $T$ as claimed.

Let $\phi : T\longrightarrow W$ be the restriction of $\phi_p$
to $T$. Since there is an irreducible component $S$ of $\phi^{-1}(V)$
with $\dim\,S\geq(\dim\,T-\dim\,W)+\dim\,V$, we deduce
$$\mld(\eta_V;X,Y)=\mld(V;X,Y)\leq\tau+\codim(V,W)$$
 via another application of Theorem~\ref{char_log_canonical}.
This concludes the proof.
\end{proof}

%References
\providecommand{\bysame}{\leavevmode \hbox \o3em
{\hrulefill}\thinspace}

\end{document}